\documentclass{dcds-B}
\usepackage{graphics}

\def\currentvolume{1}
\def\currentissue{1}
\def\currentyear{2001}
\def\currentmonth{January}
\def\ppages{1--14}

\def\evec{{\mbox{\boldmath{$e$}}}}
\def\hvec{{\mbox{\boldmath{$h$}}}}
\def\mvec{{\mbox{\boldmath{$m$}}}}
\def\uvec{{\mbox{\boldmath{$u$}}}}
\def\vvec{{\mbox{\boldmath{$v$}}}}
\def\wvec{{\mbox{\boldmath{$w$}}}}
\def\Hvec{{\mbox{\boldmath{$H$}}}}
\def\Mvec{{\mbox{\boldmath{$M$}}}}

\title {Hysteresis in layered spring magnets}
\author [J. Samuel Jiang,
         Hans G. Kaper,
         and
         Gary K. Leaf]{}


\email{jiang@anl.gov}
\email{kaper@mcs.anl.gov}
\email{leaf@mcs.anl.gov}

\begin{document}

\maketitle

\centerline{\scshape J.~Samuel Jiang}
\medskip

{\footnotesize
   \centerline{ Materials Science Division}
   \centerline{ Argonne National Laboratory}
   \centerline{ Argonne, IL 60439} }
\medskip

\centerline{\scshape Hans G.~Kaper} 
\medskip

{\footnotesize
   \centerline{ Mathematics and Computer Science Division }
   \centerline{ Argonne National Laboratory}
   \centerline{ Argonne, IL 60439} }
\medskip

\centerline{\scshape Gary K.~Leaf} 
\medskip

{\footnotesize
   \centerline{ Mathematics and Computer Science Division }
   \centerline{ Argonne National Laboratory}
   \centerline{ Argonne, IL 60439} }
\medskip

\centerline{(Communicated by Shouhong Wang)}

\medskip

\begin{abstract}
This article addresses a problem of micromagnetics:
the reversal of magnetic moments in layered spring magnets.
A one-dimensional model is used of a film consisting of
several atomic layers of a soft material on top of
several atomic layers of a hard material.
Each atomic layer is taken to be uniformly magnetized,
and spatial inhomogeneities within an atomic layer are neglected.
The state of such a system is described by a chain
of magnetic spin vectors.
Each spin vector behaves like a spinning top
driven locally by the effective magnetic field
and subject to damping (Landau--Lifshitz--Gilbert equation).
A numerical integration scheme for the LLG equation
is presented that is unconditionally stable and preserves
the magnitude of the magnetization vector at all times.
The results of numerical investigations for a bilayer
in a rotating in-plane magnetic field show
hysteresis with a basic period of $2\pi$
at moderate fields and
hysteresis with a basic period of $\pi$
at strong fields.
\end{abstract}

\section{Introduction\label{s-intro}} 
Exchange-spring coupled magnets
(\textit{spring magnets}, for short)
hold significant promise for applications
in information recording and storage devices.
Spring magnets consist of nanodispersed hard and soft
magnetic phases that are coupled at the interfaces.
(In a hard material, the magnetic moment tends
to be aligned with the easy axis;
in a soft material, it is more or less free
to align itself with the local magnetic field.)
The superior magnetic properties of a spring magnet
stem from the fact that the soft phase enhances
the magnetization of the composite~\cite{Kneller-Hawig-91,
Coey-Skomski-93,Skomski-Coey-93,Fischer-L-K-98,
Fullerton-et-al-98,Jiang-et-al-00}.
Since the performance of a spring magnet is determined
by the stability of the soft phase against magnetization reversal,
it is important to identify the factors affecting
the reversal process.

Thin films provide an interesting class of simple models
for which one can perform both physical and computational
experiments.
A spring-magnet structure can be realized by interleaving
hard and soft magnetic layers, and because the magnetic
properties vary predominantly along the normal direction,
the structure of such spring magnets is essentially
one dimensional.

In this article we investigate magnetic reversal in a
hard/soft bilayer---a layer of soft material on top of
a layer of hard material---with strong coupling at the interface.
The hard and soft layers both consist of several atomic layers;
each atomic layer is treated as uniformly magnetized,
and spatial inhomogeneities within an atomic layer are neglected.
The state of the bilayer is thus described by a chain of spins,
each spin representing the magnetic moment of an atomic layer.

The dynamics of a magnetic moment are entirely local.
A magnetic moment is like a spinning top, which is driven
by the effective magnetic field and subject to damping.
The relevant equation was first formulated by
Landau and Lifshitz~\cite{LL35} and later given
in an equivalent form by Gilbert~\cite{G55}.
The local effective field is derived variationally
from an energy functional~\cite{FS00}.

The Landau--Lifshitz--Gilbert equation
preserves the magnitude of the magnetic moment,
and it is important to maintain this property
in numerical approximations.
Our first purpose in this article is to present
an integration scheme that conserves magnetization
at all times.

Our second purpose is to present some results
of numerical simulations for a bilayer.
They show two types of rotational hysteresis:
one at moderate fields with a basic period of $2\pi$,
which is associated with the irreversible behavior
of the chirality (``handedness'') of the chain of
magnetic spins in the soft layers;
another at strong fields with a basic period of $\pi$,
which is associated with the irreversible behavior
of a full-length transition of the chain of magnetic spins
in both the hard and the soft layers.
These results explain the experimental observation
of hysteresis in some torque measurements~\cite{torque}.
They also agree qualitatively with some magneto-optical
measurements of the magnetization angle~\cite{MMM}.
However, they differ at the quantitative level.
The discrepancy is due to the one-dimensional model,
which does not allow for the nucleation and motion
of nanodomains.

The remainder of this article consists of three sections.
In Section~\ref{s-model}, we describe the mathematical model,
together with the numerical approximation procedure.
In Section~\ref{s-results}, we present some simulation
results for hard/soft bilayers.
(Additional results are presented in~\cite{JKL}.)
In Section~\ref{s-summary}, we summarize our conclusions.

\section{Mathematical model\label{s-model}}
A layered spring magnet is a multilayer structure,
which consists of $N_h$ atomic layers of a hard magnetic material
adjacent to $N_s$ atomic layers of a soft magnetic material,
\begin{eqnarray*}
\mbox{Hard layers}:&& i \in I_h = \{1, \ldots\,, N_h \} , \\
\mbox{Soft layers}:&& i \in I_s = \{N_h + 1, \ldots\,, N_h + N_s \} .
\end{eqnarray*}
We put $I = I_h \cup I_s$ and $N = N_h + N_s$.
The atomic layers are homogeneous, and variations occur only
in the direction normal to the layers.
We assume for convenience that the atomic layers are equally thick;
their thickness $d$ is of the order of angstroms
(1~\AA~$ = 1.0 \cdot 10^{-8}$~cm).

We adopt a right-handed Cartesian $(x, y, z)$ coordinate system,
where the $x$ and $y$ axes are in the plane of an atomic layer,
the $x$ axis coincides with the easy axis of the hard material,
and the $z$ axis is in the direction normal to the layers;
$\evec_x$, $\evec_y$, and $\evec_z$ are the unit vectors
in the direction of increasing $x$, $y$, and $z$, respectively.
In a polar $(\phi, \theta)$ coordinate system,
$\phi$ is the out-of-plane angle and $\theta$ the
in-plane angle measured counterclockwise
from the positive $x$ axis.

The state of the bilayer is completely described
by the set of \textit{magnetic moments},
\begin{equation}
  \Mvec = \{ \Mvec_i : i \in I \} .
  \label{M}
\end{equation}
Each $\Mvec_i$ is a vector-valued function of time $t$,
with components $M_{i,x}$, $M_{i,y}$, and $M_{i,z}$.
The magnitude $M_i$ of $\Mvec_i$ is the
\textit{magnetization},
the unit vector $\mvec_i = \Mvec_i/M_i$ is
the \textit{magnetic spin} in the $i$th layer.
The magnetization is constant at all times
and equal to the local saturation magnetization,
\begin{equation}
  \Mvec_i (t)
  =
  M_i \mvec_i (t) ,
  \quad\mbox{with }
  M_i
  =
  \left\{
  \begin{array}{ll}
  M_h & \mbox{ if } i \in I_h , \\
  M_s & \mbox{ if } i \in I_s .
  \end{array}
  \right.
  \label{Mi}
\end{equation}
Here, $M_h$ and $M_s$ are the values of
the saturation magnetization for the
hard and soft material, respectively.
Each magnetic spin can be specified in terms
of its Cartesian or polar components,
\begin{equation}
  \mvec_i
  =
  (m_{i,x}, m_{i,y}, m_{i,z})^{\scriptsize{\mbox{t}}}
  =
  (\cos \phi_i \cos \theta_i,
   \cos \phi_i \sin \theta_i,
   \sin \phi_i)^{\scriptsize{\mbox{t}}} .
  \label{mi}
\end{equation}
Thus, $\theta_i$ is the \textit{in-plane} angle
of $\mvec_i$ with the easy axis of the hard material
(measured from the positive $x$ direction),
$\phi_i$ the \textit{out-of-plane} angle of $\mvec_i$.

\subsection{Dynamics of the magnetic moment\label{ss-dynamics}}
A magnetic moment is like a spinning top,
which is driven by the local effective magnetic field
and subject to damping.
The equation of motion for $\Mvec_i$ is the
Landau--Lifshitz--Gilbert (LLG) equation,
\begin{equation}
  \frac{\partial \Mvec_i}{\partial t}
  = - \gamma (\Mvec_i \times \Hvec_i)
  + \frac{g}{M_i} \left(\Mvec_i \times \frac{\partial \Mvec_i}{\partial t}\right) ,
  \quad i \in I .
  \label{LLG-M}
\end{equation}
Here, $\Hvec_i$ is the effective magnetic field
in the $i$th layer,
$\gamma$ the gyromagnetic constant, and
$g$ a (dimensionless) damping coefficient.
Note that the LLG~equation yields a magnetic moment
whose magnitude is constant in time.
An equivalent form of the LLG~equation is
\begin{equation}
  \frac{\partial \Mvec_i}{\partial t}
  = - c \left[
  (\Mvec_i \times \Hvec_i) + \frac{g}{M_i} \Mvec_i \times (\Mvec_i \times \Hvec_i)
  \right] ,
  \quad i \in I ,
  \label{LLG}
\end{equation}
where $c = \gamma / (1 + g^2)$.
We rescale $t$ by a factor $c$ and take $c = 1$ from here on.

The effective magnetic field is found by taking
the variational derivative of the free energy,
\begin{equation}
  \Hvec_i
  = \Hvec_a
  - \frac{\delta F}{\delta \Mvec_i} ,
\end{equation}
where $\Hvec_a$ is the externally applied field,
which we take to be uniform and constant in time.
The free-energy density $F$ is the sum of
the exchange energy,
the anisotropy energy,
and the demagnetization energy,
\begin{equation}
  \mathcal{F} [\Mvec]
  =
  \int_\Omega
  \left[
  \frac{1}{2} A(z) \left| \frac{\partial \mvec}{\partial z} \right|^2
  +
  K(z) \left| \mvec \times \evec_x \right|^2
  + \frac{1}{2}(4\pi) ( \Mvec \cdot \evec_z)^2
  \right] .
  \label{F}
\end{equation}
Here, $\Omega$ is the $z$ interval occupied
by the entire multilayered structure,
$A$ is the exchange coupling coefficient, and
$K$ is the anisotropy coefficient.
The demagnetization tensor for a layer
has only one element, $D_{zz}$;
$4\pi$ is its value for an infinitely thin
flat ellipsoid~\cite{Bertram}.
In practice, one approximates $\Hvec_i$
by the expression
\begin{eqnarray}
  \Hvec_i
  =
  \Hvec_a
  &\hspace{-0.7em}+\hspace{-0.7em}&
  \frac{1}{M_i}
  \left[
  J_{i, i+1} (\mvec_{i+1} - \mvec_i)
  - J_{i, i-1} (\mvec_i - \mvec_{i-1})
  \right]
  - 2 \frac{K_i}{M_i} \evec_x \times (\mvec_i \times \evec_x)
  \nonumber \\
  &\hspace{-0.7em}-\hspace{-0.7em}& 4\pi M_i (\mvec_i \cdot \evec_z) \evec_z ,
  \quad i \in I ,
  \label{Hi}
\end{eqnarray}
where
\begin{equation}
  \mvec_0 = \mvec_1 , \; \mvec_{N+1} = \mvec_N .
\end{equation}
The coupling coefficient $J$,
which is related to $A$ ($J = Ad^{-2}$),
has the same value between layers of the same material;
similarly, the anisotropy coefficient $K$
is constant within the same material,
\begin{equation}
  J_{i, i+1} = \left\{
  \begin{array}{ll}
  J_h ,    & i = 1, \ldots\,, N_h - 1 , \\
  J_{hs} , & i = N_h , \\
  J_s ,    & i = N_h + 1, \ldots\,, N ,
  \end{array}
  \right. \;
  K_i = \left\{
  \begin{array}{ll}
  K_h , & i = 1, \ldots\,, N_h , \\
  K_s , & i = N_h + 1, \ldots\,, N .
  \end{array}
  \right.
\end{equation}
The actual values of these material parameters
depend on the temperature;
$K_s \ll K_h$ in all practical cases.

\subsection{Integration of the LLG equation\label{ss-integration}}
The LLG equation maintains a constant magnetization,
so the only quantity that changes in the course of time
is the direction of the magnetic moment.
We therefore begin by rewriting the LLG equation
in terms of $\mvec$.
As the equation is entirely local to each layer,
we drop the index $i$ temporarily.
We use the prime ${}'$ to denote differentiation
with respect to time.

Let $H$ be the strength of the magnetic field,
and let $\hvec = \Hvec / H$ be the unit vector
in the direction of $\Hvec$,
\begin{equation}
  \Hvec (t) = H (t) \hvec (t) .
\end{equation}
Then the LLG equation is
\begin{equation}
  \mvec'
  = - H
  \left[
  (\mvec \times \hvec) + g \mvec \times (\mvec \times \hvec)
  \right] .
  \label{LLG-m}
\end{equation}
We decompose the equation by means of
the projection operators $P$ and $Q$,
\begin{equation}
  P\uvec = (\uvec \cdot \hvec) \hvec , \;
  Q\uvec = \uvec - P\uvec = \hvec \times (\uvec \times \hvec) ,
  \quad \uvec \in {\bf R}^3 .
\end{equation}
Equation~(\ref{LLG-m}) is equivalent to
the two equations
\begin{eqnarray}
  P\mvec'
  &\hspace{-0.7em}=\hspace{-0.7em}& - H
  P \left[
  (\mvec \times \hvec) + g \mvec \times (\mvec \times \hvec)
  \right] ,
  \label{Pm'} \\
  Q\mvec'
  &\hspace{-0.7em}=\hspace{-0.7em}& -H
  Q \left[
  (\mvec \times \hvec) + g \mvec \times (\mvec \times \hvec)
  \right] .
  \label{Qm'}
\end{eqnarray}
Notice the identities
\begin{eqnarray}
  &P (\mvec \times \hvec) = \mathbf{0} , \;
  P [\mvec \times (\mvec \times \hvec)] = (\mvec \cdot Q\mvec) \hvec
  = - [1 - (P\mvec \cdot P\mvec)^2 ] \hvec ,&\mbox{} \\
  &Q (\mvec \times \hvec) = - J Q\mvec , \;
  Q [\mvec \times (\mvec \times \hvec)] = (\mvec \cdot \hvec) Q\mvec ,&\mbox{}
\end{eqnarray}
where $J$ is the square root of the negative identity in $\mathbf{R}^2$,
\begin{equation}
  I = \left( \begin{array}{rr} 1&0 \\ 0&1 \end{array} \right) , \;
  J = \left( \begin{array}{rr} 0&-1 \\ 1&0 \end{array} \right) , \;
  J^2 = - I .
\end{equation}
Hence, we can recast Eqs.~(\ref{Pm'}) and~(\ref{Qm'})
in the form
\begin{eqnarray}
  P\mvec'
  &\hspace{-0.7em}=\hspace{-0.7em}& gH [1 - (P\mvec \cdot P\mvec)^2 ] \hvec ,
  \label{Pm'red} \\
  Q\mvec'
  &\hspace{-0.7em}=\hspace{-0.7em}& H
  [J - g (\mvec \cdot \hvec) I ]  Q\mvec .
  \label{Qm'red}
\end{eqnarray}
Suppose that the direction of $\Hvec$ does not change
on an interval $(t, t + \Delta t)$,
\begin{equation}
  \hvec(s) = \hvec(t) , \quad s \in (t, t + \Delta t) .
\end{equation}
Then
$P\mvec' = (P\mvec)'$ and $Q\mvec' = (Q\mvec)'$
on $(t, t + \Delta t)$,
so Eqs.~(\ref{Pm'red}) and (\ref{Qm'red})
reduce to a coupled system of differential equations
for the scalar $u = (P \mvec \cdot \hvec)$ in $\mathbf{R}$
and the vector $\vvec = Q \mvec$ in $\mathbf{R}^2$,
\begin{eqnarray}
  u' &\hspace{-0.7em}=\hspace{-0.7em}& gH (1 - u^2)
  \quad \mbox{on } (t, t + \Delta t) ,
  \label{u'} \\
  \vvec' &\hspace{-0.7em}=\hspace{-0.7em}& H ( J - g u I ) \vvec
  \quad \mbox{on } (t, t + \Delta t) .
  \label{v'}
\end{eqnarray}
From these equations we conclude
that the critical states are
$u = 1$, $\vvec = \mathbf{0}$
($\mvec = \hvec$, magnetic moment parallel
to the magnetic field)
and $u = -1$, $\vvec = \mathbf{0}$
($\mvec = - \hvec$, magnetic moment
antiparallel to the magnetic field).
The former is linearly stable, the latter unstable
under infinitesimal perturbations.

We now turn to the integration of Eqs.~(\ref{u'}) and~(\ref{v'}).
The former is independent of $\vvec$
and can be integrated immediately.
If not only the direction, but also the magnitude of $H$
is constant on $(t, t + \Delta t)$,
\begin{equation}
  \Hvec(s) = \Hvec(t) ,
  \quad s \in (t, t + \Delta t) ,
\end{equation}
we find
\begin{equation}
  u(s)
  =
  \frac{u(t) \cosh(gH(t)(s - t)) + \sinh(gH(t)(s - t))}
       {\cosh(gH(t)(s - t)) + u(t) \sinh(gH(t)(s - t))} ,
  \quad s \in (t, t + \Delta t) .
  \label{u}
\end{equation}
Next, we turn to Eq.~(\ref{v'}).
We replace the constant $gH$ by $u'/(1-u^2)$ (from Eq.~(\ref{u'}))
and use the identity $-uu'/(1-u^2) = (\ln (1-u^2)^{1/2})'$
to convert the equation into a differential equation
for the vector $\wvec = (1 - u^2)^{-1/2} \vvec$,
\begin{equation}
  \wvec' = H J \wvec \quad \mbox{on } (t, t + \Delta t) .
  \label{w'}
\end{equation}
This equation can be integrated,
\begin{eqnarray}
  \wvec(s)
  &\hspace{-0.7em}=\hspace{-0.7em}&
   \mathrm{e}^{H(t) (s - t) J} \wvec(t) \nonumber \\
  &\hspace{-0.7em}=\hspace{-0.7em}&
   [ \cos(H(t)(s - t)) I + \sin(H(t) (s - t)) J ] \wvec (t) ,
  \quad s \in (t, t + \Delta t) .
\end{eqnarray}
From the expression~(\ref{u}) we obtain
\begin{equation}
  (1-u(s)^2)^{1/2}
  =
  \frac{(1 - u(t)^2)^{1/2}}
       {\cosh(gH(t) (s - t)) + u(t) \sinh(gH(t)(s - t))} ,
\end{equation}
so
\begin{equation}
  \vvec(s)
  =
  \frac{\cos(H(t)(s - t)) I + \sin (H(t) (s - t)) J}
       {\cosh (gH(t) (s - t)) + u(t) \sinh (gH(t) (s - t))}
  \vvec(t) ,
  \quad s \in (t, t + \Delta t) .
  \label{v}
\end{equation}
These results suggest the following choice of the integration scheme
for Eq.~(\ref{LLG-m}):
\begin{eqnarray}
  \mvec_{n+1}
  &\hspace{-0.7em}=\hspace{-0.7em}&
  \frac{(\mvec_n \cdot \hvec_n) \cosh(gH_n \Delta t)
        + \sinh(gH_n \Delta t)}
       {\cosh(gH_n \Delta t) + (\mvec_n \cdot \hvec_n) \sinh(gH_n \Delta t)}
  \hvec_n \nonumber \\
  &\hspace{-1.0em}&\mbox{}+
  \frac{\cos(H_n \Delta t) I + \sin (H_n \Delta t) J}
       {\cosh(gH_n \Delta t) + (\mvec_n \cdot \hvec_n) \sinh(gH_n \Delta t)}
  \hvec_n \times (\mvec_n \times \hvec_n) ,
  \label{m+}
\end{eqnarray}
where
$\mvec_{n+1} = \mvec(t_{n+1})$,
$\mvec_n = \mvec(t_n)$,
$\hvec_n = \hvec(t_n)$,
$H_n = H(t_n)$, and
$\Delta t = t_{n+1} - t_n$.

The algorithm~(\ref{m+}) is unconditionally stable
for all values of $\Delta t$.
Of course, the quality of the approximation suffers
as $\Delta t$ increases.
However, the algorithm explicitly displays
the relationship between the size of $\Delta t$
and the local error in the time integration.
The rate of precession of $\mvec$ around
the polar axis is governed by $H$, the
magnitude of the local effective field:
in one time step, $\mvec$ precesses
through an angle $H \Delta t$.
Therefore, by properly choosing $\Delta t$,
we can resolve the fastest precessional motion
in a given number of time steps per period.
Since $H$ varies over the course of a simulation,
we have a natural and direct means to adjust the size
of $\Delta t$ to the current dynamical state,
while maintaining the resolution of the precessional motion.

Other algorithms for the numerical integration
of the LLG equation have been proposed recently by
Nigam~\cite{nigam} and E and Wang~\cite{E}.

\subsection{Computing equilibrium configurations\label{ss-equil}}
The analysis in the preceding section
suggests the following algorithm for finding
equilibrium spin configurations.
Starting from a given equilibrium state
$\Mvec = \{ \Mvec_i : i \in I \}$
at time $t_0$,
one uses Eq.~(\ref{Hi}) to compute the magnetic field
$\Hvec_i$ in each layer at $t_0$.
Having found $\Hvec_i (t_0)$ for all $i \in I$,
one advances in time to $t_1 = t_0 + \Delta t$
and uses Eqs.~(\ref{Mi}) and (\ref{m+})
to compute $\Mvec$ at $t_1$.
If $\Delta t$ is sufficiently small,
$\Mvec (t_1)$ is a close approximation of
the state of the system at time $t_1$.
One continues this process, finding approximations
at successive times
$t_n = t_0 + n \Delta t$, $n = 1, 2, \ldots$\,,
until equilibrium is reached.

\section{Numerical results\label{s-results}}
The algorithm of the preceding section has been used
to study hysteresis phenomena in hard/soft bilayers
that are driven by an applied field~$\Hvec_a$
that is uniform, constant in time, and parallel
to the planes of the atomic layers.
The expression for the effective magnetic field,
Eq.~(\ref{Hi}), decomposes into an in-plane component,
\begin{eqnarray}
  \Hvec_i \times \evec_z
  =
  \Hvec_a \times \evec_z
  &\hspace{-0.7em}+\hspace{-0.7em}&
  \frac{1}{M_i}
  \left[
  J_{i, i+1} (\mvec_{i+1} - \mvec_i)
  - J_{i, i-1} (\mvec_i - \mvec_{i-1})
  \right] \times \evec_z
  \nonumber \\
  &\hspace{-0.7em}-\hspace{-0.7em}& 2 \frac{K_i}{M_i} (\mvec_i \cdot \evec_y) \evec_x ,
  \quad i \in I ,
  \label{Hi-in}
\end{eqnarray}
and an out-of-plane component,
\begin{eqnarray}
  \Hvec_i \cdot \evec_z
  &\hspace{-0.7em}=\hspace{-0.7em}&
  \frac{1}{M_i}
  \left[
  J_{i, i+1} (\mvec_{i+1} - \mvec_i)
  - J_{i, i-1} (\mvec_i - \mvec_{i-1})
  \right]
  \cdot \evec_z
  \nonumber \\
  &&\mbox{}- 2 \frac{K_i}{M_i} \mvec_i
  \cdot \evec_z
  - 4\pi M_i \mvec_i 
  \cdot \evec_z ,
  \quad i \in I .
  \label{Hi-out}
\end{eqnarray}
When the system is in an equilibrium state,
the effective magnetic field is parallel
(or antiparallel) to the magnetic spin;
see Section~\ref{ss-integration}.
Hence, each $\Hvec_i$ is a multiple of $\mvec_i$,
and Eq.~(\ref{Hi-out}) reduces to a homogeneous system
of linear algebraic equations for the set of scalars
$\{\mvec_i \cdot \evec_z: i \in I\}$.
In general, this system admits only the trivial solution,
so the magnetic moments lie in the plane of the atomic layers.
In the notation of Eq.~(\ref{mi}),
$\phi_i = 0$ for all $i \in I$ at equilibrium,
and the only relevant variables are the in-plane angles
$\{ \theta_i : i \in I\}$.
(Of course, the magnetic spin may have an out-of-plane
component during the transient phase of the computation.)

In the numerical simulations we focus on the
in-plane angle of the magnetic spin at equilibrium
and investigate its behavior as a function of
the strength $H_a$ and the direction $\theta_a$
of the applied field,
\begin{equation}
  \Hvec_a = H_a \hvec_a , \quad
  \hvec_a = (\cos \theta_a, \sin \theta_a, 0)^{\scriptsize{\mbox{t}}} .
\end{equation}
The following computations refer to a bilayer configuration
consisting of
$N_h = 115$ atomic layers of Sm-Co (a hard material) and
$N_s = 100$ atomic layers of Fe (a soft material).
A different configuration is used in Section~\ref{ss-compare},
where we make a comparison with some magneto-optical measurements.
Table~\ref{t-parms} gives the values of
the material parameters $A$, $K$, and $M$,
as well as the values of the coupling coefficient $J$
($J = Ad^{-2}$, $d = 2$~\AA).
In all cases, the damping coefficient $g = 0.5$.
\begin{table}[h]
\caption{Numerical values of the parameters.}  \label{t-parms}
\vspace{-3ex}
\begin{footnotesize}
\begin{center}
\begin{tabular}{||l||c|c|c|c||}\hline
                & $A$ (erg/cm)        & $J$ (erg/cm$^3$)   & $K$ (erg/cm$^3$) & $M$ (emu/cm$^3$) \\\hline
Fe              & $2.8 \cdot 10^{-6}$ & $7.0 \cdot 10^{9}$ & $1.0 \cdot 10^3$ & 1,700            \\
Interface       & $1.8 \cdot 10^{-6}$ & $4.5 \cdot 10^{9}$ & --               & --               \\
Sm-Co           & $1.2 \cdot 10^{-6}$ & $3.0 \cdot 10^{9}$ & $5.0 \cdot 10^7$ & 550              \\\hline
\end{tabular}
\end{center}
\end{footnotesize}
\end{table}

\subsection{Rotational hysteresis\label{ss-fixed}}
The case $H_a = 4800$~oersteds is typical,
at least for moderate values of $H_a$
(see Section~\ref{ss-variable}).

\begin{center}
\begin{figure}[htb]
\resizebox{2.3in}{!}{\mbox{\includegraphics{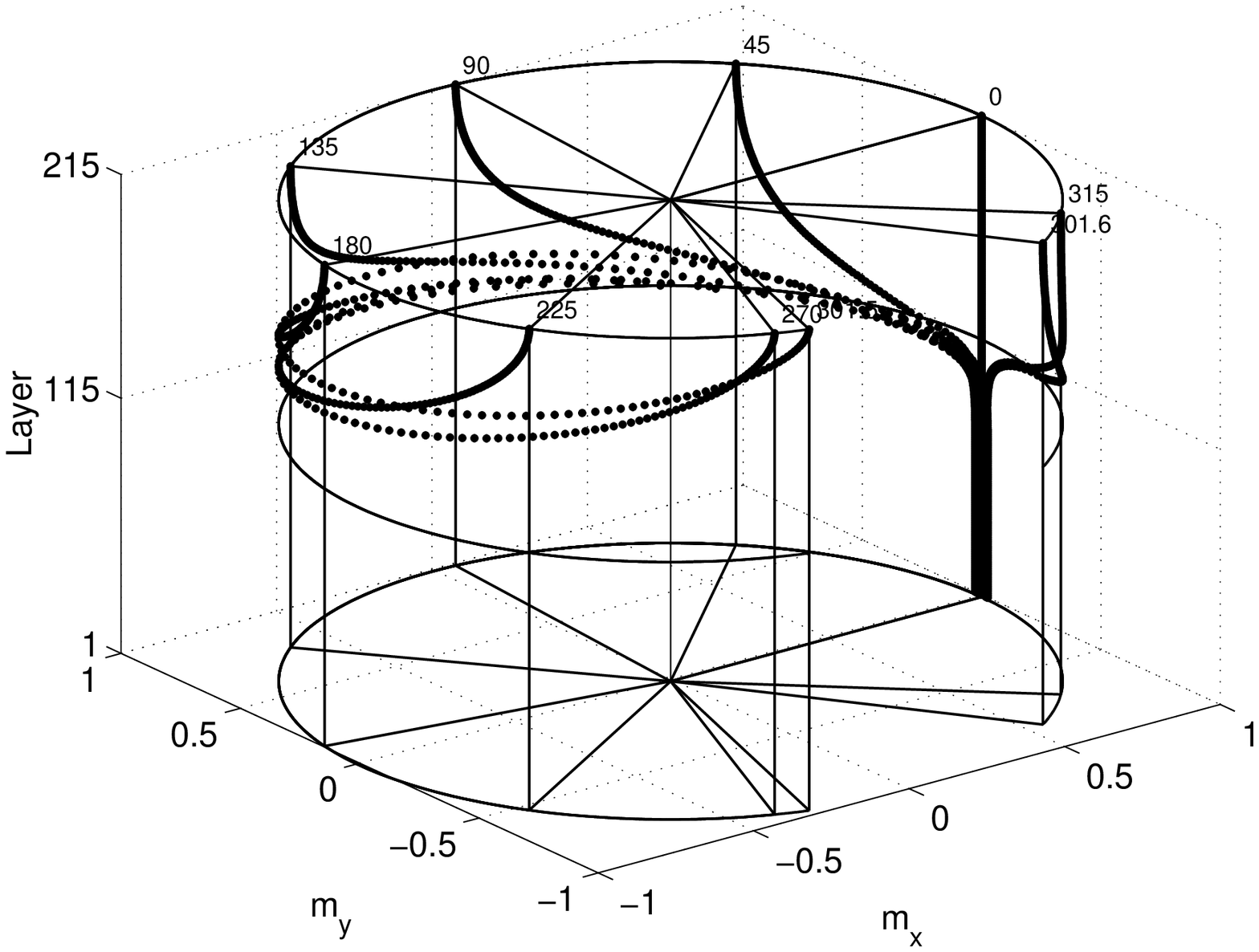}}}
\hfill
\resizebox{2.3in}{!}{\mbox{\includegraphics{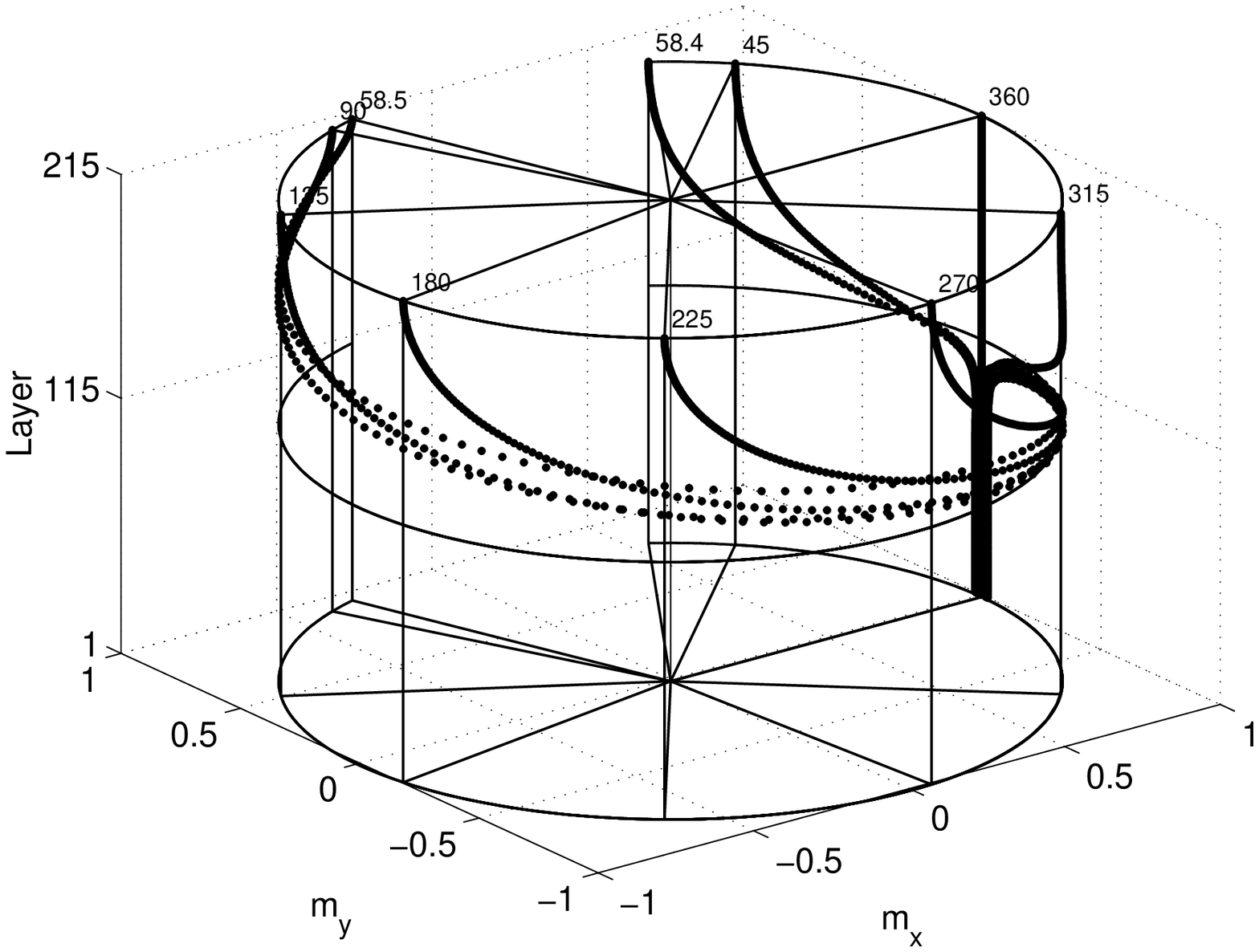}}}
\vspace{-2ex}
\caption{Equilibrium spin configurations; $H_a = 4800$~oersteds.
         Left: $\theta_a$ increasing, right: $\theta_a$ decreasing.
  \label{f-spin}}
\end{figure}
\end{center}

The simulations show that
the equilibrium spin configurations
for increasing $\theta_a$ ($0 < \theta_a < 2\pi$)
and decreasing $\theta_a$ ($2\pi > \theta_a > 0$)
are mirror images of each other.
Figure~\ref{f-spin} shows two sets of magnetic spin configurations
at equilibrium for various values of $\theta_a$,
one set (left) as $\theta_a$ increases from 0 to $2\pi$,
the other set (right) as $\theta_a$ decreases from $2\pi$ to 0.
The heavy dots represent the endpoints of the magnetic spin
(a unit vector) in each layer for various angles $\theta_a$;
the values of $\theta_a$, in degrees, are indicated near the top layer.
(The dots merge into a solid line where the magnetic spins
in adjacent layers are close.)

\begin{center}
\begin{figure}[htb]
\center\resizebox{2.4in}{!}{\mbox{\includegraphics{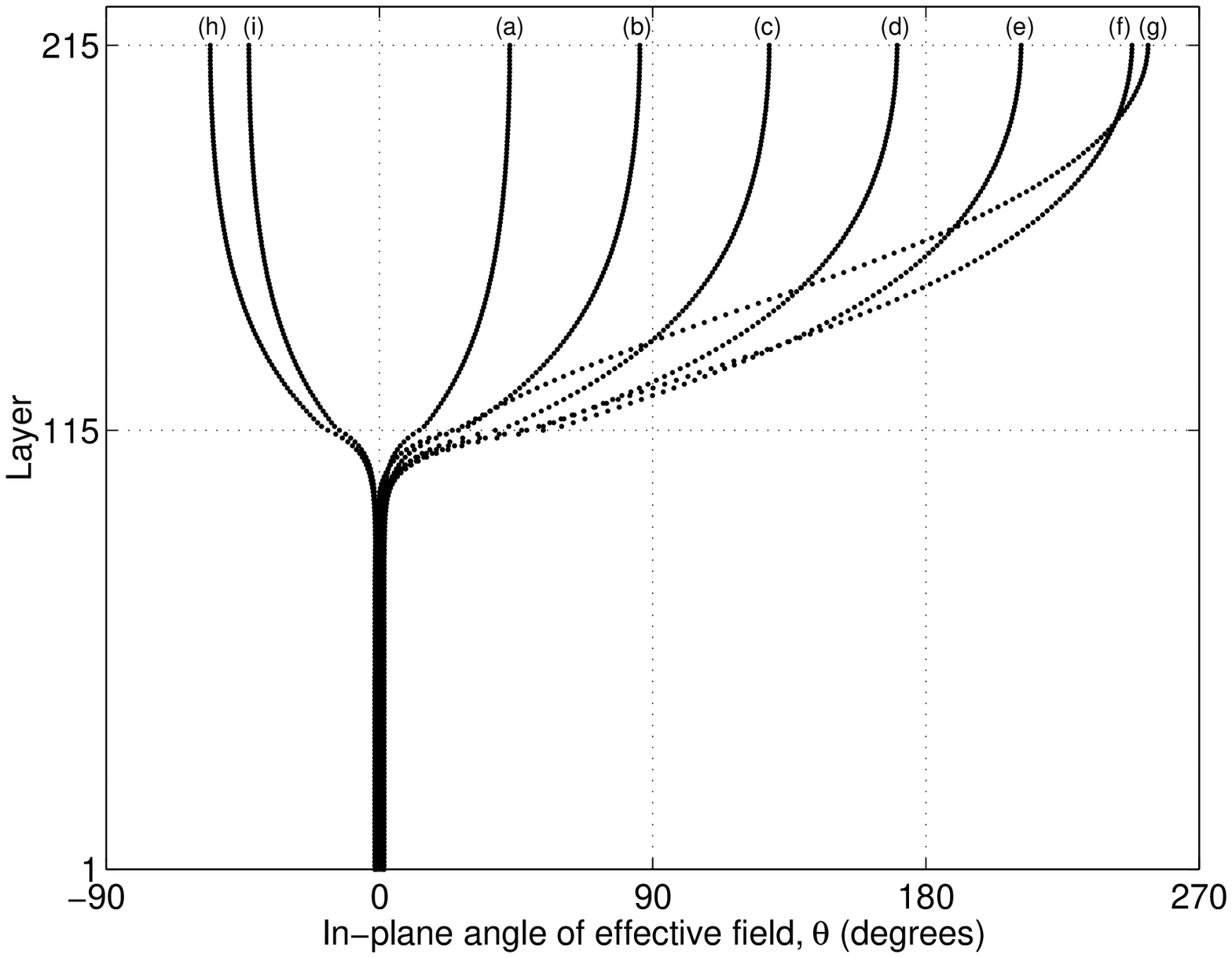}}}
\vspace{-2ex}
\caption{In-plane angle $\theta_i$ vs.~$i$; $H_a = 4800$~oersteds;
         (a)~$\theta_a = 45$,
         (b)~$\theta_a = 90$,
         (c)~$\theta_a = 135$,
         (d)~$\theta_a = 180$,
         (e)~$\theta_a = 225$,
         (f)~$\theta_a = 270$,
         (g)~$\theta_a = 301.5$,
         (h)~$\theta_a = 301.6$,
         (i)~$\theta_a = 315$~degrees.
  \label{f-theta_z}}
\end{figure}
\end{center}

Notice that the chirality (``handedness'')
of the chain of magnetic spins changes
from positive at $\theta_a = 301.5$
to negative at $\theta_a = 301.6$~degrees
and from negative at $\theta_a = 58.5$
to positive at $\theta_a = 58.4$~degrees.
Figure~\ref{f-theta_z} shows this change
in a different way.

Here, we have plotted the in-plane angle
$\theta_i$ against the layer index $i$
for increasing values of $\theta_a$.
(The graphs for decreasing values of
$\theta_a$ are obtained by symmetry.)
First, the graph changes continuously (but not
monotonically) as $\theta_a$ increases from 0
to 301.5~degrees, $\theta_i$ increasing
with $i$ (positive chirality).
Then it changes discontinously
as $\theta_a$ increases to 301.6~degrees:
$\theta_i$ suddenly becomes decreasing
instead of increasing with $i$.
Finally, it changes continuously again
as $\theta_a$ increases further,
$\theta_i$ decreasing with $i$
(negative chirality),
to return to the original graph
($\theta_i = 0$ for all $i \in I$)
as $\theta_a$ reaches 360~degrees.
In all cases, the spin is fixed
along the easy axis ($\theta_i = 0$)
in most of the hard layers;
it begins to deviate from the easy axis
only as one approaches the interface ($i = 115$).
The first derivative is discontinuous
at the interface, and the tangent is vertical 
in the top layer ($i = 215$).

\begin{center}
\begin{figure}[htb]
\center\resizebox{2.4in}{!}{\mbox{\includegraphics{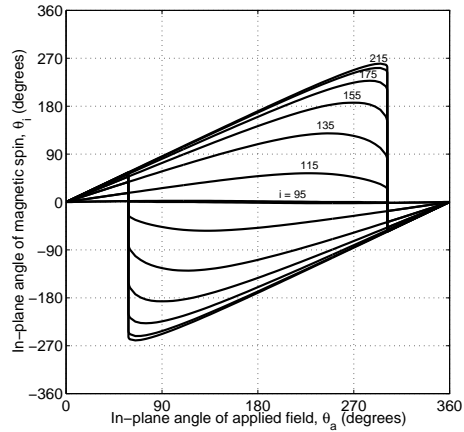}}}
\vspace{-2ex}
\caption{Rotational hysteresis:
         in-plane angle $\theta_i$ vs.\ $\theta_a$;
         $H_a = 4800$~oersteds;
         $i = 95, 115, 135, 155, 175, 195, 215$.
  \label{f-hyst}}
\end{figure}
\end{center}

The change in chirality is irreversible and
induces \textit{rotational hysteresis}.
The in-plane angle of each spin vector
traverses a different trajectory
as the applied field rotates 360~degrees
in the forward and backward direction.
The hysteresis loop has the same shape,
and particularly the same width, in all layers.
Its vertical dimension contracts gradually
as one descends through the soft layers,
to disappear entirely in the hard layers
somewhat below the interface;
see Fig.~\ref{f-hyst}.

\subsection{Two types of rotational hysteresis\label{ss-variable}}
\begin{center}
\begin{figure}[htb]
\center\resizebox{2.4in}{!}{\mbox{\includegraphics{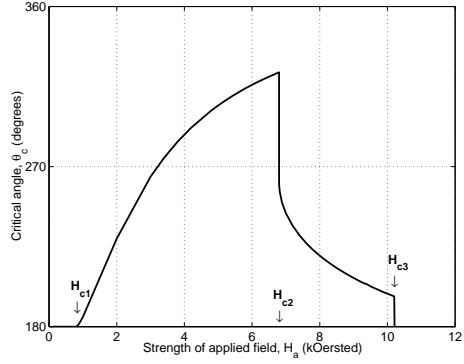}}}
\vspace{-2ex}
\caption{The critical angle $\theta_c$ as a function of $H_a$.
  \label{f-thetac}}
\end{figure}
\end{center}

When the strength of the applied field is varied,
we observe different modes of irreversible behavior.
We recall (Fig.~\ref{f-hyst})
that, as $\theta_a$ increases from 0,
the chirality changes discontinuously from
positive to negative as the direction of the applied field
deviates sufficiently from the easy axis.
We denote the critical value of the angle $\theta_a$
by $\theta_c$
($\theta_c = 301.5\ldots$ at $H_a = 4800$~oersteds).
Figure~\ref{f-thetac} shows the variation of
$\theta_c$ with $H_a$.

As long as $H_a$ is sufficiently small,
the magnetization process is reversible.
At a first critical value of $H_a$,
marked $H_{c1}$,
the chirality of the chain of magnetic spins
changes for the first time, and
rotational hysteresis of the type
discussed in the preceding section
(with a basic period of 360 degrees)
sets in.
The width of the hysteresis loop,
which is symmetric around $\theta_a = \pi$,
increases monotonically from 0 at $H_a = H_{c1}$
to some value less than $2\pi$.

\begin{center}
\begin{figure}[htb]
\center\resizebox{2.4in}{!}{\mbox{\includegraphics{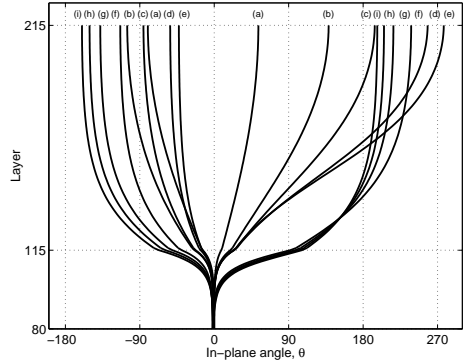}}}
\vspace{-2ex}
\caption{In-plane angle $\theta_i$ vs.~$i$;
  (a)~$H_a = 1000$, (b)~$H_a = 2000$, (c)~$H_a = 3000$,
  (d)~$H_a = 5000$, (e)~$H_a = 6000$, (f)~$H_a = 7000$,
  (g)~$H_a = 8000$, (h)~$H_a = 9000$, (i)~$H_a = 10,000$~oersteds.
  Right branches: $\theta_a$ just below~$\theta_c$,
  left branches: $\theta_a$ just above~$\theta_c$.
  \label{f-tulip}}
\end{figure}
\end{center}

At a second critical value of $H_a$,
marked $H_{c2}$,
a sharp discontinuity occurs.
The hysteresis loop narrows significantly and
continues to narrow as $H_a$ increases beyond $H_{c2}$.
The cause of this discontinuity becomes obvious
in Fig.~\ref{f-tulip}, where we have plotted
$\theta_i$ against~$i$;
cf.~Fig.~\ref{f-theta_z}.
(The bottom 80 layers of hard material, where $\theta_i$ does not
deviate noticeably from 0, are not included in this figure.)
As $H_a$ reaches the value $H_{c2}$,
the chain of spins has been stretched to its widest extent;
it can no longer support the span in the top layer,
stiffens suddenly, and becomes more like a rigid rod.
The rod-like behavior is apparent from the increasing range
where the chain is almost vertical.

\begin{center}
\begin{figure}[htb]
\center\resizebox{2.4in}{!}{\mbox{\includegraphics{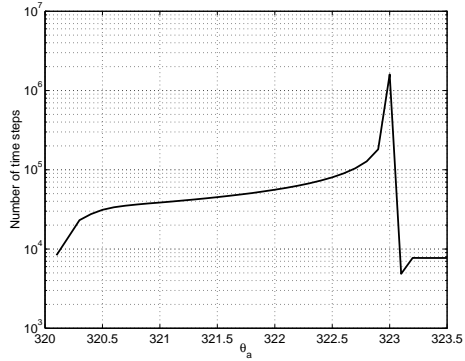}}}
\vspace{-2ex}
\caption{Equilibration time near $\theta_c$;
  $H_a = 6797$~oersteds, $\theta_c = 323.0$~degrees.
  \label{f-eq_time}}
\end{figure}
\end{center}

The structural change in the chain of spins
has some of the characteristics of a phase transition.
For example, we observe a significant increase in the equilibration time
(by two orders of magnitude) as $\theta_a$ approaches $\theta_c$;
see Fig.~\ref{f-eq_time}.
Also, the increasing size of the rigid domain near $H_{c2}$
is reminiscent of a diverging correlation length.

\begin{center}
\begin{figure}[htb]
\center
\resizebox{2.4in}{!}{\mbox{\includegraphics{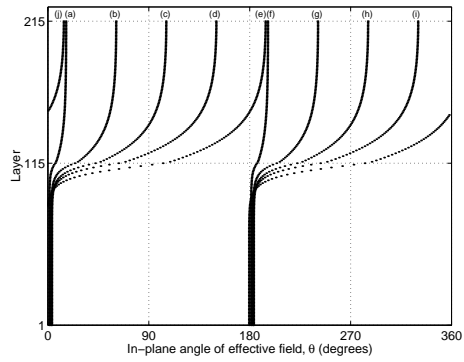}}}
\vspace{-2ex}
\caption{In-plane angle $\theta_i$ vs.~$i$; $H_a = 10,400$~oersteds;
         (a)~$\theta_a = 16$,
         (b)~$\theta_a = 61$,
         (c)~$\theta_a = 106$,
         (d)~$\theta_a = 151$,
         (e)~$\theta_a = 195$,
         (f)~$\theta_a = 196$,
         (g)~$\theta_a = 241$,
         (h)~$\theta_a = 286$,
         (i)~$\theta_a = 331$,
         (j)~$\theta_a = 375$~degrees.
  \label{f-snap}}
\end{figure}
\end{center}

At a third critical value of $H_a$,
marked $H_{c3}$,
another significant change occurs.
The field has now become sufficiently strong
to move the spins in both the soft and the hard materials.
From here on, the chain of spin vectors changes
over its entire length, maintaining its chirality.
Figure~\ref{f-snap} gives $\theta_i$ vs.\ $i$
for increasing values of $\theta_a$.
(The graphs for decreasing values of $\theta_a$
are obtained by symmetry.)
The value $H_a = 10,400$~oersteds is just above $H_{c3}$.
This figure should be compared with Fig.~\ref{f-theta_z}
for the standard case, $H_a = 4800$~oersteds.
The exact determination of $H_{c3}$ is delicate;
in our numerical simulations we found a slight
rate dependence in the regime near $H_{c3}$.

\begin{center}
\begin{figure}[htb]
\center
\resizebox{2.4in}{!}{\mbox{\includegraphics{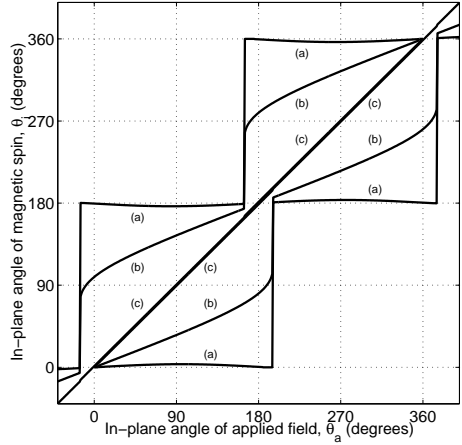}}}
\vspace{-2ex}
\caption{In-plane angle $\theta_i$ vs.~$\theta_a$; $H_a = 10,400$~oersteds;
         (a)~$i = 85$,
         (b)~$i = 115$,
         (c)~$i = 215$.
  \label{f-upandup}}
\end{figure}
\end{center}

Because the chain of spins behaves more
like an elastic spring than a stiff rod,
a new type of rotational hysteresis emerges,
whose basic period can be any multiple of 180 degrees.
Figure~\ref{f-upandup} shows three graphs:
one graph (c) is along the diagonal;
the other two (a and b) are symmetric
with respect to the diagonal.
The outer graph (a) shows $\theta_i$
for $i = 85$ (hard layer).
The part below the diagonal is traversed
in the upward direction
as $\theta_a$ increases from 0;
the part above the diagonal is traversed
in the downward direction
as $\theta_a$ decreases from 360~degrees.
The spin is oriented in either the positive
or the negative $x$ direction.
Transitions occur at $\theta_c$ and
at every multiple of $\pi$ beyond $\theta_c$.
The center graph (c) shows
$\theta_i$ for $i = 215$ (top layer).
The orientation of this spin varies
continuously with $\theta_a$ and
is perfectly reversible.
Finally, the middle graph (b) shows
$\theta_i$ for $i = 115$ (at the interface).
Here, the spin rotates continuously until it jumps.
The jumps occur at $\theta_c$ and
at every multiple of $\pi$ beyond $\theta_c$.
The graphs for the remaining layers fill the space
between the ones drawn in the figure.
The main point to observe is that the graphs for
$\theta_a$ increasing always increase and stay
below the diagonal, while those for $\theta_a$ decreasing
always decrease and stay above the diagonal.
Hence, chirality is preserved in both cases.

When the direction of $\theta_a$ is reversed,
$\theta_i$ crosses the diagonal as soon as $\theta_a - \theta_c$
is a multiple of $\pi$; after crossing, it remains
on the part of the graph situated on the newly reached
side of the diagonal.
Because there is a gap between the graphs for $\theta_i$
in the interior layers and the diagonal,
the orientation of the magnetic spin shows
rotational hysteresis in all interior layers.
This hysteresis is caused by a full-length transition
of the chain of magnetic spins, rather than the partial-range
transition that was responsible for the hysteresis
below $H_{c3}$.

\subsection{Comparison with experiment\label{ss-compare}}
\begin{center}
\begin{figure}[htb]
\resizebox{2.37in}{!}{\mbox{\includegraphics{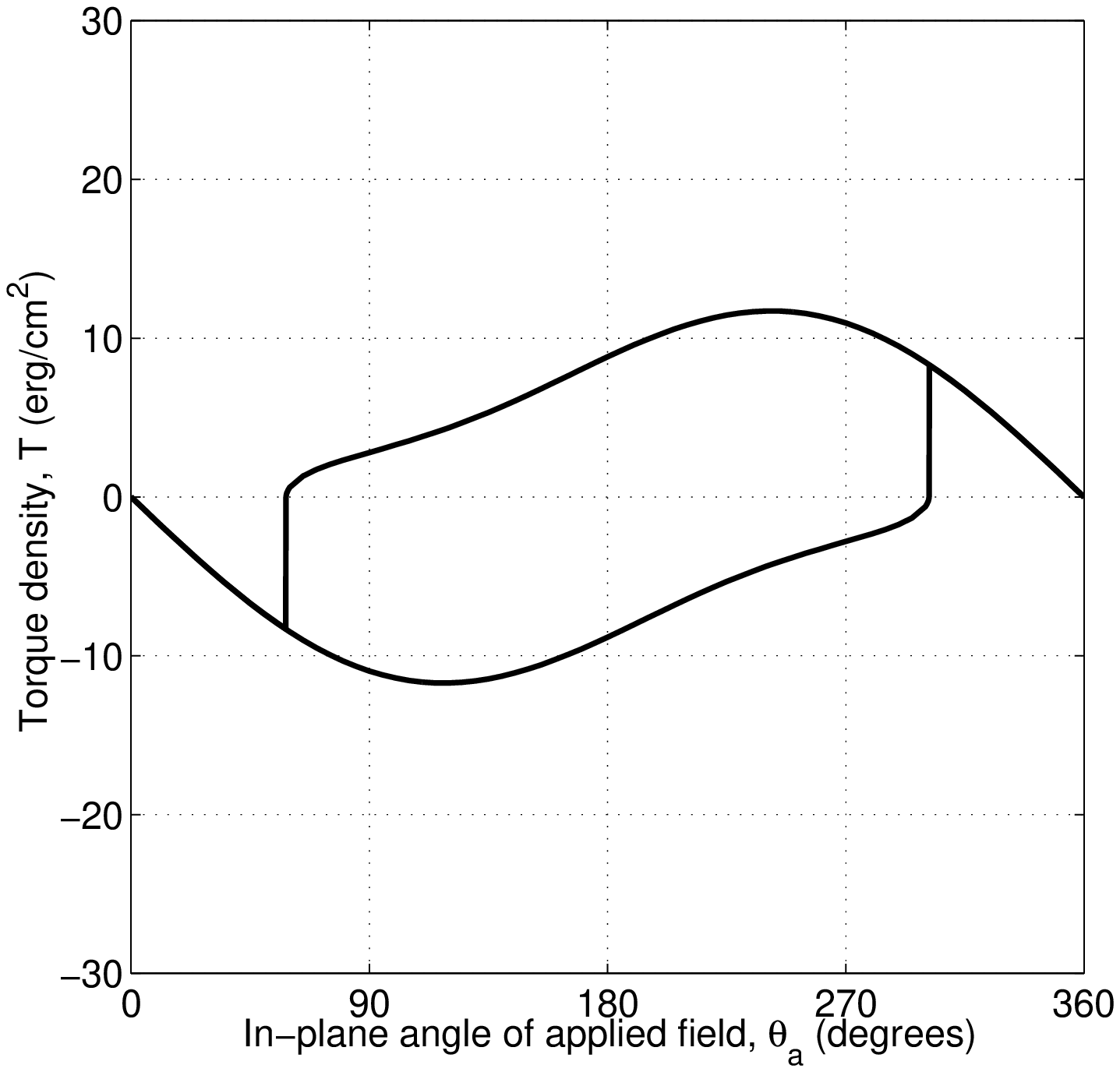}}}
\hfill
\resizebox{2.3in}{!}{\mbox{\includegraphics{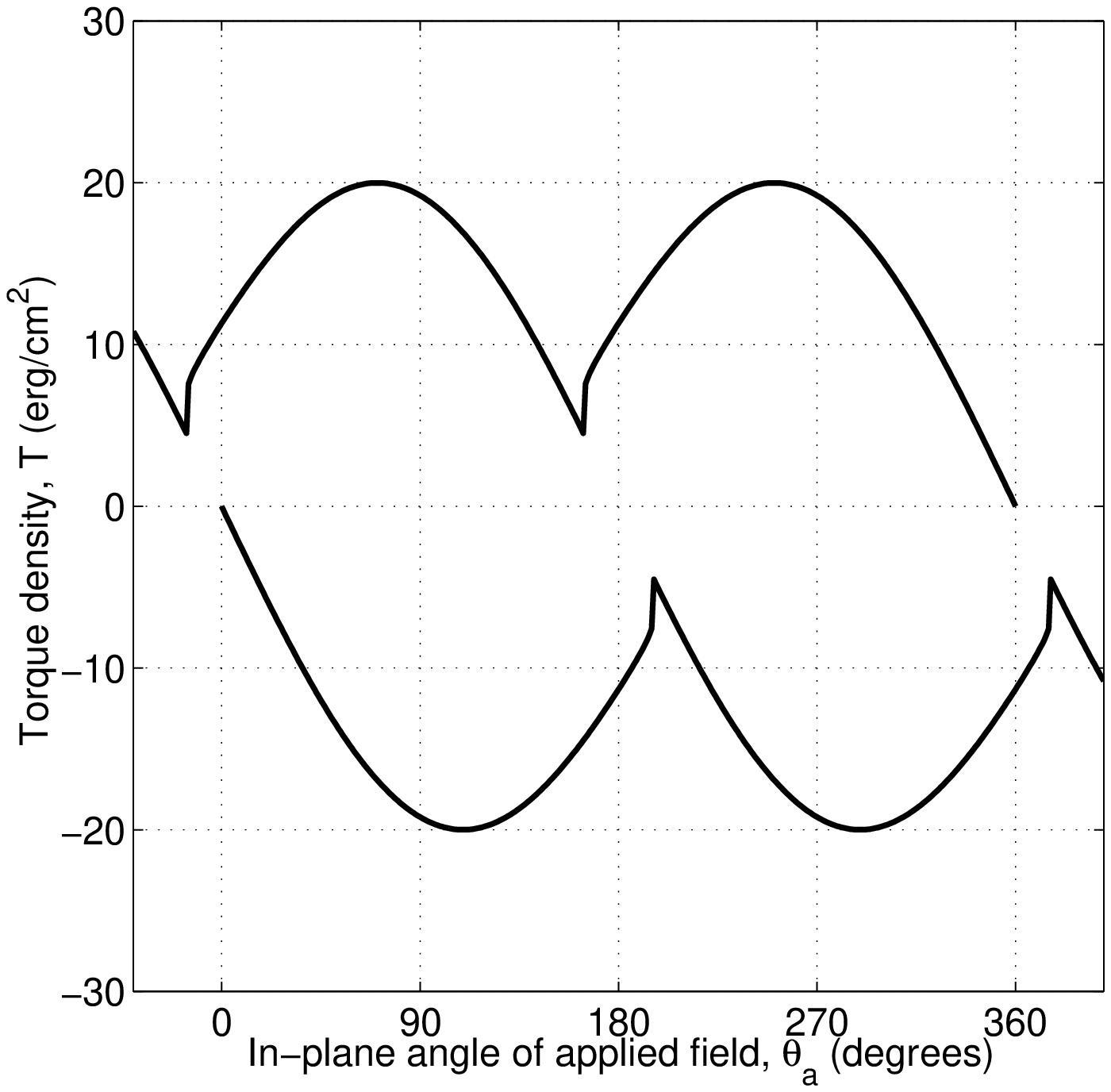}}}
\vspace{-2ex}
\caption{Torque density; $H_a = 4800$ (left) and $10,400$ (right) oersteds.
  \label{f-torque}}
\end{figure}
\end{center}

Quantities such as the magnetic moment are fundamental
to describe the state of the system, but they are
not directly measurable in an experiment.
Measurable quantities are 
the \textit{torque density} $T$ and
the \textit{magnetization angle} $\alpha$,
\begin{equation}
  T = H_a d \sum_{i \in I} M_i \sin (\theta_a - \theta_i) ,
  \quad
  \alpha
  =
  \tan^{-1} \frac {\sum_{i \in I} M_i \sin \theta_i}
                  {\sum_{i \in I} M_i \cos \theta_i} .
  \label{alpha}
\end{equation}
Both $T$ and $\alpha$ reflect the hysteretic behavior
of the magnetic moments.
Figure~\ref{f-torque} shows the torque density computed at
$H_a = 4800$ and $10,400$~oersteds.

Experimental torque measurements at comparable values
of $H_a$ show similarly shaped graphs,
with extrema at approximately the same values of $\theta_a$,
but significantly narrower hysteresis loops~\cite{torque}.

\begin{center}
\begin{figure}[htb]
\resizebox{2.4in}{!}{\mbox{\includegraphics{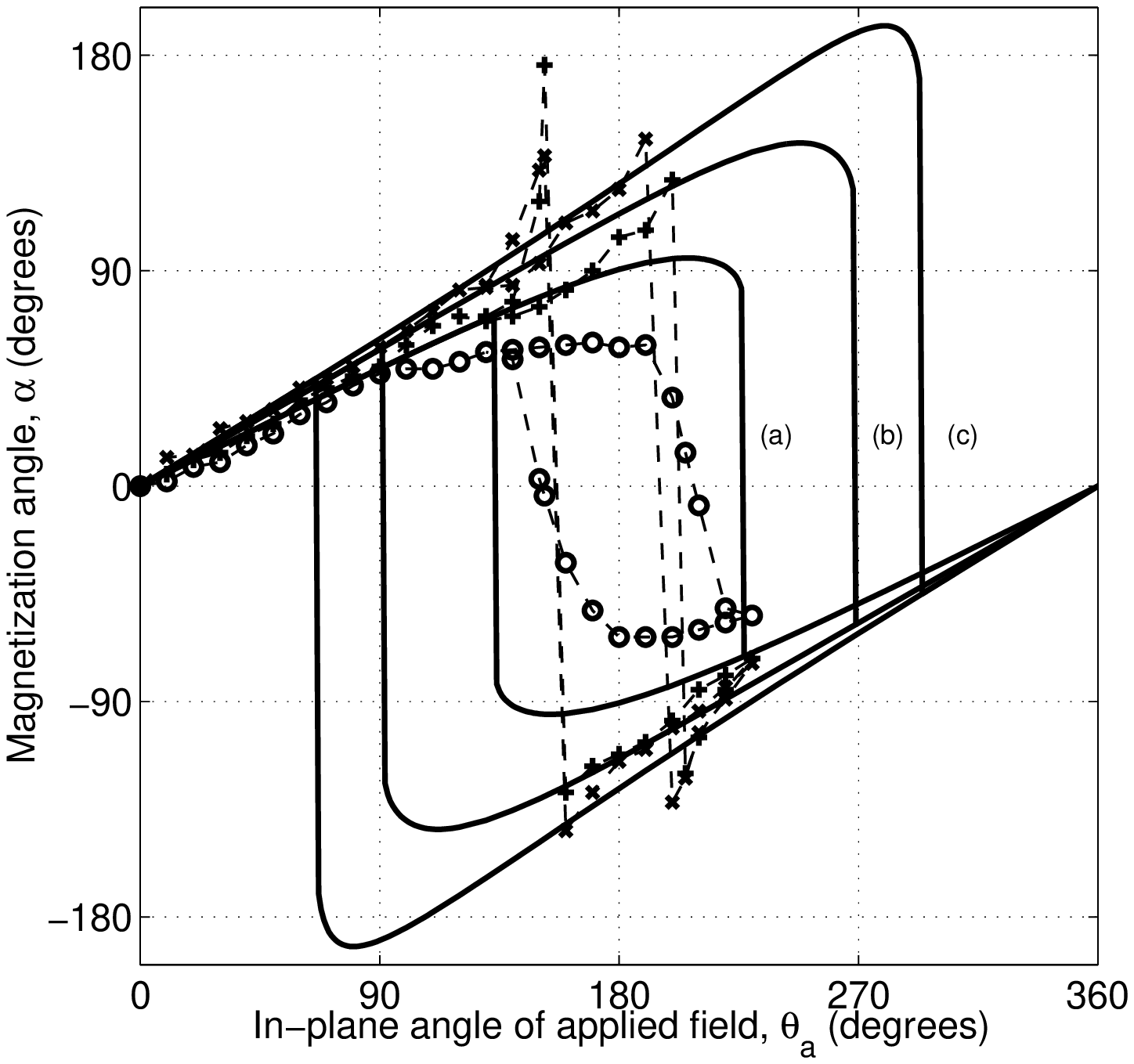}}}
\vspace{-2ex}
\caption{Magnetization angle; (a)~$H_a = 360$ (\textbf{o}),
                              (b)~$H_a = 600$ ($\mathbf{+}$), and
                              (c)~$H_a = 840$ ($\mathbf{\times}$)~oersteds.
  \label{f-angle}}
\end{figure}
\end{center}

In Fig.~\ref{f-angle}, we compare results for
the magnetization angle with experimental data.
The data were obtained by magneto-optical means
for a bilayer consisting of
$N_h = 100$ atomic layers of Sm-Co and
$N_s = 250$ atomic layers of Fe;
the simulation curves also refer to this
configuration~\cite{MMM}.
The measurements were done at relatively low fields
($H_a = 360, 600$, and 840~oersteds)
and for a limited range of directions
($\theta_a = 0:10:230$~degrees).

There is certainly qualitative agreement,
but the simulations generally yield
wider hysteresis loops than
the experiments, and
the discrepancy becomes greater
as the field strength increases.
This behavior can be explained by the fact that
the model used in the simulations
is a single-domain model,
which does not allow for the important phenomenon
of nucleation and motion of nanodomains.
As a result, the demagnetization energy is
seriously overestimated.
In realistic simulations,
one must use multidimensional models and
allow for lateral inhomogeneities~\cite{MMM}.

\section{Conclusions\label{s-summary}}
In this article we have addressed an important issue
in micromagnetics: magnetization reversal in
layered spring magnets.
We have used a one-dimensional model of a film
consisting of atomic layers of a soft material
on top of atomic layers of a hard material
with strong coupling at the interface,
assuming no variation in the lateral directions.
The state of such a system is described by
a chain of magnetic spin vectors.
Each spin vector behaves like a spinning top
driven by the local magnetic field and subject
to damping.
The dynamics are described by a system of
LLG equations, Eq.~(\ref{LLG}), coupled with
a variational equation for the magnetic field,
Eq.~(\ref{Hi}).

We have presented an integration procedure
that maintains the invariance of the magnetization
(the magnitude of the magnetization vector)
and proposed an algorithm for finding the
equilibrium state of the system.

We have applied the algorithm to simulate
magnetization reversal in layered spring magnets.
The results show that a layered spring magnet exhibits
rotational hysteresis with a basic period of 360~degrees
at moderately strong fields and rotational hysteresis
with a basic period of 180~degrees at strong fields.
The former type of hysteresis is induced
by a partial-length transition of the chain
of magnetic spins;
the transition occurs only in the soft material
and causes a change of chirality.
The hysteresis in strong fields is induced by
a full-length transition of the chain of spins
in both the hard and the soft layers;
it is much weaker than the rotational hysteresis
at moderately strong fields and can cover any period
that is a multiple of the basic period.

The numerical results for the torque and
magnetization angle agree qualitatively
with the experimental data but differ
at the quantitative level.
In particular, the one-dimensional model seriously
overestimates the demagnetization energy, since
it does not allow for the nucleation and motion
of nanodomains.
In realistic simulations, lateral inhomogeneities
must be taken into account.

\section*{Acknowledgments}
This work was supported by the Mathematical, Information, and
Computational Sciences Division subprogram of the Office of Advanced
Scientific Computing Research, U.S.~Department of Energy, under
Contract W-31-109-Eng-38.
Most of the numerical simulations were carried out by
Jaime Hernandez Jr.\ (University of Texas at El Paso),
who was a participant in the Energy Research Undergraduate
Laboratory Fellowship program at Argonne National Laboratory
(summer 2000).

\smallskip

Received January 2001.

\end{document}